\documentclass[review,sort&compress,colorlinks=true]{elsarticle}

\usepackage{natbib}

\usepackage{amssymb}
\usepackage{epsfig}
\usepackage{graphicx}
\usepackage{amsmath,amsthm,bm,mathrsfs}
\usepackage{epsf}

 \newtheoremstyle{theorem}{6pt}{6pt}{\rm}{}{\sffamily}{ }{ }{}
 \theoremstyle{theorem}

 \newtheoremstyle{algorithm}{6pt}{6pt}{\rm}{}{\sffamily}{ }{ }{}
 \theoremstyle{algorithm}

 \newtheoremstyle{lemma}{6pt}{6pt}{\rm}{}{\sffamily}{ }{ }{}
 \theoremstyle{lemma}

\newtheoremstyle{case}{6pt}{6pt}{\rm}{}{\sffamily}{. }{ }{}
 \theoremstyle{case}

 \newtheoremstyle{statement}{6pt}{6pt}{\rm}{}{\sffamily}{ }{ }{}
\theoremstyle{statement}

 \newtheoremstyle{corollary}{6pt}{6pt}{\rm}{}{\sffamily}{ }{ }{}
 \theoremstyle{corollary}

  \newtheoremstyle{definition}{6pt}{6pt}{\rm}{}{\sffamily}{ }{ }{}
 \theoremstyle{definition}

\newtheoremstyle{example}{6pt}{6pt}{\rm}{}{\sffamily}{ }{ }{}
\theoremstyle{example}

\newtheoremstyle{remark}{6pt}{6pt}{\rm}{}{\sffamily}{ }{ }{}
\theoremstyle{remark}

\newtheoremstyle{approximation}{6pt}{6pt}{\rm}{}{\sffamily}{ }{ }{}
\theoremstyle{approximation}

\newtheoremstyle{scheme}{6pt}{6pt}{\rm}{}{\sffamily}{ }{ }{}
\theoremstyle{scheme}

\newtheoremstyle{Algorithm}{6pt}{6pt}{\rm}{}{\sffamily}{ }{ }{}
\theoremstyle{Algorithm}

\newtheoremstyle{Assumption}{6pt}{6pt}{\rm}{}{\sffamily}{ }{ }{}
\theoremstyle{Assumption}

\newtheoremstyle{proposition}{6pt}{6pt}{\rm}{}{\sffamily}{ }{ }{}
\theoremstyle{proposition}

\newtheoremstyle{hypo}{6pt}{6pt}{\rm}{}{\sffamily}{ }{ }{}
 \theoremstyle{hypo}

  \newtheoremstyle{Step}{6pt}{6pt}{\rm}{}{}{ }{ }{}
 \theoremstyle{Step}

\numberwithin{equation}{section}

\begin{document}

\title{Performance enhancement of non-minimum phase feedback systems by fractional-order cancellation of non-minimum phase zero on the Riemann surface: New theoretical and experimental results}
\author{ {\sc Farshad Merrikh-Bayat and Aliakbar Salimi}\\[2pt]
Department of Electrical and Computer Engineering, University of Zanjan, \\
 Zanjan, IRAN, Email: $\{\mathrm{f.bayat,a\_salimi\}@znu.ac.ir}$.\\[6pt]
\vspace*{6pt}}
\pagestyle{headings}
\markboth{F. Merrikh-Bayat and A. Salimi}{\rm PERFORMANCE ENHANCEMENT OF NON-MINIMUM PHASE FEEDBACK SYSTEMS}
\maketitle

\address{Department of Electrical and Computer Engineering, University of Zanjan, Zanjan, Iran}


\begin{abstract}
{The non-minimum phase (NMP) zero of a linear process located in
the feedback connection cannot be cancelled by the same pole of
controller according to the internal instability problem. However,
such a zero can partly be cancelled by the same fractional-order pole of a pre-compensator located in series with process without facing internal instability. This paper
first presents new theoretical results on the properties of this method of cancellation, and provides design techniques for the pre-compensator. It is especially
shown that by appropriate design of pre-compensator this method can simultaneously increase the gain and
phase margin of the system under control without a considerable reduction of open-loop bandwidth, and consequently, it can make the control problem easier to
solve. Then, a method for realization of such a pre-compensator is proposed and performance of the resulted closed-loop system is studied through an experimental setup.
}
{non-minimum phase zero, performance limitation, partial zero-pole cancellation, experimental result, Riemann surface.}
\end{abstract}

\section{Introduction}
In the field of linear time-invariant systems, a process is
identified as a non-minimum phase (NMP) system if its transfer
function has at least one right half-plane (RHP) zero, or a RHP
pole, or a time delay. Among others, systems with RHP zero(s)
constitute a very important category of NMP systems, both from
theoretical and practical point of view. Such zeros appear in many
real-world systems such as flexible link robots \citep{vakil,
benosman}, step-up DC-DC converters \citep{zhang, chen},
floating-wind turbines \citep{Fischer}, aircrafts \citep{chemori},
bicycles \citep{astrom1}, driving a car backwards \citep{hoagg},
continuous stirred tank reactors \citep{tofighi}, nano positioning
devices \citep{aggarwal}, and many others.

One classical fact in relation to NMP zeros is that they cannot be
cancelled by the same poles of controller according to the
internal instability problem \citep{skogestad}. The other
well-known fact is that NMP zeros of the process put some
limitations on the performance of the corresponding feedback
system \citep{sidi}. One main reason for this limitation is that in
order to achieve a good command following and disturbance
rejection behavior, any feedback system needs large open-loop
gains at lower frequencies, which is often provided by the
controller \citep{skogestad}. On the other hand, according to the
classical root-locus method, any closed-loop system has the
property that its poles move towards open-loop zeros as the gain
in the loop is increased. Hence, there is a tradeoff between
performance and stability when the process has a NMP zero since
the open-loop gain cannot be increased arbitrarily in this case.
It should be noted that many other controller design techniques
also strictly depend on the possibility of applying large gains in
the loop. For example, successful application of the loop transfer
recovery (LTR) method needs using large gains in the loop
\citep{skogestad}. That is why application of LTR is mainly limited to
minimum phase processes.

The difficulty of controlling NMP processes can also be explained
through frequency domain analysis. More precisely, a process with
a NMP zero is more phase-lag compared to the one which has a zero
with the same amplitude but at the left half-plane. This extra
phase-lag puts a limitation on the gain of controller since
increasing the loop gain increases the gain crossover frequency
which often leads to decreasing the PM. Again, considering the
fact that large open-loop gains at low frequencies are essential
for command following and disturbance rejection,
it is obvious that NMP zero puts a serious restriction on the
performance of the closed-loop system.

Very recently the idea of partial cancellation of the NMP zero of
process on the Riemann surface is proposed in \citep{farshad1}. It
is especially shown in that paper that partial cancellation of the
NMP zero of process leads to a transfer function with
fractional-order zero which can be controlled more easily compared
to the original NMP process. However, no routine design technique is presented in that paper. This paper completes the basic idea
proposed in \citep{farshad1} both from theoretical and experimental
aspects. The main idea behind the methods proposed in this paper
for partial cancellation of NMP zero is that in dealing with many
real-world processes, smaller the value of
$\kappa\triangleq|P(j\omega_u)/P(0)|$ more easier its control
\citep{astrom2}, where $P(s)$ is the transfer function of process
and $\omega_u$ is the ultimate frequency, {\it i.e.} $\angle
P(j\omega_u)=-180^\circ$ (as a rule of thumb, processes with
$\kappa<0.4$ are considered as easier control problems
\citep{astrom2}). Hence, the aim of this paper is to design a
pre-compensator for partial cancellation of NMP zero such that the
resulted system has a larger gain and phase margin, and
consequently, be easier to control.

The rest of this paper is organized as the following. Some
preliminaries, which will be instrumental in the discussions of
next sections, are presented in Section 2. The method of partial
cancellation of NMP zero of process on the Riemann surface is also
briefly reviewed in this section. Some new theoretical results on
this subject are presented in Section 3. Specially, two methods
for designing the partial canceller of NMP zero are proposed in
this section. A novel NMP benchmark circuit is presented in
Section 4 and the proposed design techniques for canceller are
successfully tested on this benchmark. Finally, Section 5
concludes the paper.

\section{Preliminaries and review of previous findings}
\begin{figure}\center
\includegraphics[width=8cm]{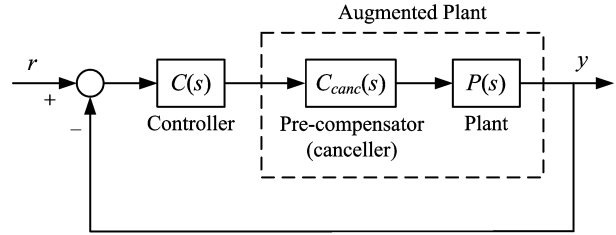}
\caption{Feedback system with pre-compensator for partial
cancellation of the NMP zero of process. $P(s)$ and $C_{canc}(s)$
are defined in (\ref{p1}) and (\ref{c_canc}),
respectively.}\label{fig_loop1}
\end{figure}

Some basic definitions, which will be instrumental in the
discussions of the next sections, are presented in this section.
Some of the material presented in this section can also be found
with more details in \citep{farshad1}. Consider the unity feedback
system shown in Fig. \ref{fig_loop1} where $P(s)$ is the NMP
transfer function of process, $C_{canc}(s)$ is the pre-compensator
used to partially cancel the NMP zero of $P(s)$ on Riemann
surface, and $C(s)$ is the controller designed to control the
\emph{augmented process} $C_{canc}(s)P(s)$. (In the rest of this
paper the terms ``canceller" and ``pre-compensator" are used to
refer to $C_{canc}(s)$ interchangeably. Moreover, in this paper
the series connection of $C_{canc}(s)$ and $P(s)$ is called the
``augmented process" since the controller $C(s)$ has to be
designed for a process with transfer function $C_{canc}(s)P(s)$
which is easier to control compared to $P(s)$ provided that
$C_{canc}(s)$ is suitably designed.) For a better and more clear
understanding the effect of proposed canceller on the function of
closed-loop system, the discussions of this paper are presented
assuming proportional controller in the loop, i.e., without a
considerable loss of generality it is assumed that the controller
$C(s)$ in Fig. \ref{fig_loop1} is in the form of a simple constant
gain. However, the main results can easily be extended to more
complicated controllers.

In the rest of this paper it is assumed that the transfer function
of NMP process in Fig. \ref{fig_loop1} can be decomposed as the
following
\begin{equation}\label{p1}
P(s)=\left(1-\frac{s}{z_{nmp}}\right) \widetilde{P}(s),
\end{equation}
where $(1-s/z_{nmp})$ is the NMP zero term of $P(s)$ and
$z_{nmp}>0$ is the NMP zero. As a very well-known classical fact
\citep{skogestad}, neither $C(s)$ nor $C_{canc}(s)$ can have a pole
at $s=z_{nmp}$ according to the internal instability problem,
i.e., any zero-pole cancellation in the RHP is impractical since
the resulted system is internally unstable. Considering the fact
that in the feedback connection of Fig. \ref{fig_loop1} the
closed-loop zeros are the same as open-loop zeros, it is concluded
that one has to tolerate the limitations caused by NMP zeros since
they appear unavoidably in the closed-loop transfer function.
However, it is shown in \citep{farshad1} that it is possible to
partly cancel the NMP zero of $P(s)$ by $C_{canc}(s)$ on the
Riemann surface and arrive at a higher performance feedback system. More precisely,
similar to \citep{farshad1}, here the transfer function of
canceller is considered as
\begin{equation}\label{c_canc}
C_{canc}(s)=\frac{1}{\sum_{k=1}^n\left(\frac{s}{z_{nmp}}\right)^{(k-1)/n}},
\end{equation}
where $n$ is the design parameter which determines the order of
cancellation of NMP zero \citep{farshad1}. (To follow the
discussions in this paper one does not need to have a deep
knowledge about fractional calculus and the time-domain
interpretation of fractional powers of $s$ in (\ref{c_canc}). See
\citep{podlubny} for more details on this subject). It is shown in
\citep{farshad1} that the series connection of $P(s)$ and
$C_{canc}(s)$, as defined in (\ref{p1}) and (\ref{c_canc}),
respectively, is as the following
\begin{equation}\label{cp}
C_{canc}(s)P(s)=\left[1-\left(\frac{s}{z_{nmp}}\right)^{1/n}\right]
\widetilde{P}(s),
\end{equation}
which, unlike $P(s)$, has a \emph{fractional-order} NMP zero at
$s=z_{nmp}$ (note that both the $P(s)$ and $C_{canc}(s)P(s)$ have
exactly the same poles and zeros, and the only difference between
these two transfer functions is that $C_{canc}(s)P(s)$ has a
fractional-order NMP zero at $s=z_{nmp}$). Assuming
\begin{equation}\label{f}
F(s)=1-\left(\frac{s}{z_{nmp}}\right)^\alpha,
\end{equation}
one can further write (\ref{cp}) as $C_{canc}(s)P(s)=F(s)
\widetilde{P}(s)$ for $\alpha= 1/n$. In brief, the effect of
canceller is that it changes the NMP term $1-s/z_{nmp}$ at the
numerator of $P(s)$ to $1-(s/z_{nmp})^\alpha$ for some rational
$\alpha$ smaller than unity.

As another definition, in the rest of this paper $L_{canc}(s)$ and
$L(s)$ denote the open-loop transfer functions with and without
applying the canceller, respectively. That is
\begin{equation}
L(s)=C(s)P(s),
\end{equation}
and
\begin{equation}
L_{canc}(s)=C(s)C_{canc}(s)P(s).
\end{equation}
The gain crossover frequency of the open-loop transfer function
with canceller is denoted as $\omega_{gc}$, that is
$|L_{canc}(j\omega_{gc})|=0$ dB.

\section{Theoretical Findings and the Proposed Design Methods for Pre-compensator}
\subsection{Effect of canceller on the open-loop frequency
response}\label{sec_pre} 
\begin{figure}\center
\includegraphics[width=9.5cm]{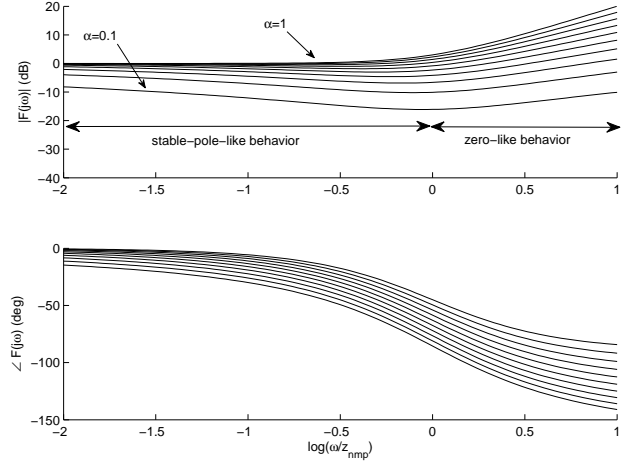}
\caption{Bode plots of the partly-cancelled zero term $F(s)=\left
[1-(s/z_{nmp})^\alpha\right]_{s=j \omega}$ ($z_{nmp}>0$) versus
the normalized frequency for $\alpha=0.1,0.2,\ldots,1$. In a
certain frequency range the magnitude and phase plots are
decreased simultaneously as the frequency is increased, which is
similar to the behavior of a stable pole.}\label{fig_freq1}
\end{figure}

One intrinsic characteristic of any zero
(either NMP or MP) is that it leads to increment in the amplitude
of the frequency response of the corresponding transfer function as the frequency is
increased. For example, in the numerator of (\ref{p1}) amplitude
of the term $\left( 1-s/z_{nmp}\right)_{s=j\omega}$ is
monotonically increased by increasing $\omega$. This increase in
the amplitude of open-loop transfer function increases the gain
crossover frequency of system and simultaneously makes it more
phase lag if the zero is NMP. Hence, it is expected that a
feedback system with a NMP zero in the loop exhibits a very poor
stability if the bandwidth is sufficiently large. But,
surprisingly, the partly-cancelled NMP term $\left
[1-(s/z_{nmp})^\alpha\right]_{s=j \omega}$ ($0<\alpha<1$) has the
property that its amplitude is monotonically \emph{decreased} by
increasing $\omega$ in the frequency range $[0,\omega_{min}]$
where $\omega_{min}= \left[\cos(\pi\alpha/2)\right] ^{1/\alpha}
z_{nmp}$ \citep{farshad1} (see Fig. \ref{fig_freq1} for more
details). It concludes that, as it can be observed in Fig.
\ref{fig_freq1}, this partly-cancelled NMP term in the numerator
of open-loop transfer function somehow acts as a \emph{stable
pole} in this frequency range. Considering the fact that in
practice the gain crossover frequency of any NMP system is often
approximately limited to the frequency of its NMP zero
\citep{horwitz, seron}, this pole-like behavior of the
partly-cancelled NMP zero can improve the gain-bandwidth of system
since it considerably decreases the amplitude of open-loop
frequency response at frequencies around $\omega=z_{nmp}$, while
it has a negligible effect at lower frequencies. In the following,
we study the frequency response of the partly cancelled zero term
(\ref{f}) mathematically assuming $\omega_{gc}=z_{nmp}$; however,
the main results can be extended to values of $\omega_{gc}\neq
z_{nmp}$ as well.

For the partly cancelled NMP term given in (\ref{f}) trivial
calculations yield
\begin{equation}\label{F1}
|F(j\omega)|_{\omega=z_{nmp}}=\left|1-\left(j\frac{\omega}
{z_{nmp}}\right)^\alpha\right|_{\omega=z_{nmp}}= \sqrt{2-2\cos
\frac{\pi\alpha}{2}}.
\end{equation}
It is concluded from (\ref{F1}) that $|F(jz_{nmp})|<1$ for
$\alpha<2/3$. It means that in the feedback system shown in Fig.
\ref{fig_loop1} application of canceller changes the NMP term of
process from $1-s/z_{nmp}$ to $1-(1/z_{nmp})^\alpha$ which,
assuming $0<\alpha<2/3$, reduces the open-loop gain at the
frequency of NMP zero (which is often close to the gain crossover
frequency of system in practice) as the following
\begin{align}
|L_{canc}(jz_{nmp})| -|L(jz_{nmp})| &=
\left|1-\frac{s}{z_{nmp}}\right|_{s=jz_{nmp}} -
\left|1-\left(\frac{s}{z_{nmp}}\right)^\alpha
\right|_{s=jz_{nmp}}\label{l_mag2}
\\&=\sqrt{2}-\sqrt{2-2\cos \frac{\pi\alpha}{2}}
\\&=\sqrt{2}\left(1-\sqrt{1-\cos\frac{\pi\alpha}{2}}\right)\label{l_mag}.
\end{align}
(Recall that, except the NMP term, $L_{canc}(s)$ is exactly the
same as $L(s)$. That is why the difference between the amplitude
of open-loop transfer functions in (\ref{l_mag2}) is expressed
only in terms of the NMP zero terms.) Reduction of open-loop gain
as given in (\ref{l_mag}) leads to reduction of gain crossover
frequency.

Two points should be noted here. First, although the amplitude of
open-loop transfer function at the frequency of NMP zero is
reduced after cancellation as given in (\ref{l_mag}), one cannot
conclude that assuming $\omega_{gc}=z_{nmp}$ the GM is also
increased in the same way. The reason is that the cancellation
also makes the system more phase lag as described below (i.e., if
$\angle L(j\omega_0)=-180^\circ$ for some $\omega_0$ then
$\angle L_{canc}(j\omega_1)=-180^\circ$ for some
$\omega_1<\omega_0$). Second, although the above discussion
studies the behavior of $|F(j\omega)|$ at $\omega=z_{nmp}$, it is
generally true that the proposed cancellation strategy reduces the
amplitude of open-loop transfer function at the gain crossover
frequency even if $\omega_{gc}\neq z_{nmp}$; however, the amount
of this reduction cannot be calculated from (\ref{l_mag}). The
proof of this statement is obvious from Fig. \ref{fig_freq1}.

Now we can study the phase behavior of the partly cancelled NMP
zero term. For the $F$ defined in (\ref{f}) we have
\begin{equation}\label{phase1}
\angle F(j\omega) \big |_{\omega=z_{nmp}} =\angle
\left[1-\left(j\frac{\omega} {z_{nmp}}\right)^\alpha
\right]_{\omega=z_{nmp}}=-\frac{\pi}{2}+\frac{\pi\alpha}{4}.
\end{equation}
Considering the fact that $\angle \left[1-s/z_{nmp} \right]
_{s=jz_{nmp}}=-\pi/4$, it is concluded from (\ref{phase1}) that
application of the fractional-order canceller makes the open-loop
system $(\alpha-1)\pi/4$ rad more phase-lag at the frequency of
NMP zero. More precisely,
\begin{equation}\label{l_phase}
\angle L_{canc}(jz_{nmp}) =\angle L(jz_{nmp})-
\frac{(1-\alpha)\pi}{4}.
\end{equation}
It is concluded from (\ref{l_mag}) and (\ref{l_phase}) that
fractional-order cancellation of NMP zero firstly decreases the
amplitude of the open-loop transfer function at the gain crossover
frequency (as well as other frequencies), which makes the system
more stable by increasing its GM, and secondly, makes the
open-loop system more phase lag at this frequency, which decreases
the PM. Hence, the proposed canceller provides us with a tradeoff
between increasing the GM and decreasing the PM. The discussions
of the next section provide us with a method to design the
canceller such that even both of the GM and PM are increased
simultaneously at the cost of slight decrement in the open-loop
gain crossover frequency (or equivalently, open-loop and closed-loop
bandwidths).

\subsection{Designing a suitable pre-compensator for partial cancellation of NMP}
Consider again the unity feedback system shown in Fig.
\ref{fig_loop1} where $P(s)$ and $C_{canc}(s)$ are defined as
given in (\ref{p1}) and (\ref{c_canc}), respectively. In the
following, two methods for designing a canceller for partial
cancellation of the NMP zero of $P(s)$ on the Riemann surface are
proposed. Before presenting these two methods it should be noted
that the main reason for using canceller is to arrive at the
augmented plant $C_{canc}(s) P(s)$ which has better properties
compared to $P(s)$ (e.g., smaller undershoot, higher PM, etc.) and
can be controlled more effectively. For this purpose, in the
following discussions only the feedback control of augmented plant
by means of a proportional controller is studied. Simplicity of
proportional controller lets us clearly understand the function of
proposed canceller. Of course, after designing the canceller and
adjusting the gain of proportional controller to the suitable
value one can use any desired method to design a new controller
for the resulted augmented plant $K_pC_{canc}(s) P(s)$, where
$K_p$ is equal to the gain of proportional controller. Hence, the
results obtained in the following can be extended to more
complicated controllers as well.

\subsubsection{Designing a canceller to increase the DC gain and keep PM unchanged}
\begin{figure}\center
\includegraphics[width=8.5cm]{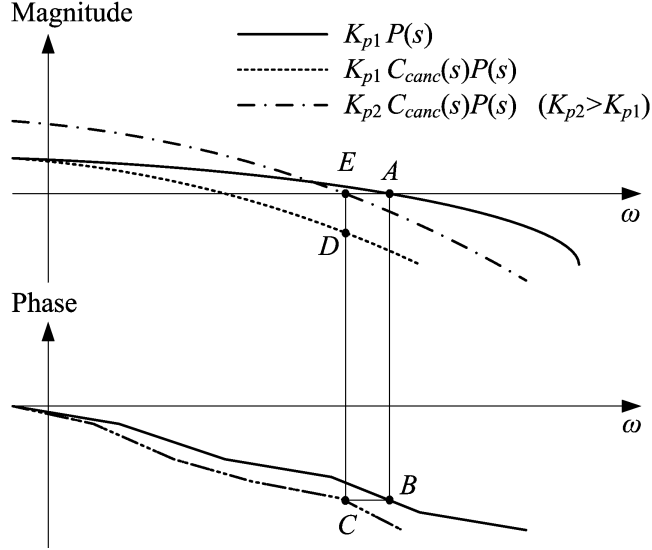}
\caption{The design procedure for canceller and proportional
controller to arrive at a feedback system with increased DC gain
(the PM remains the same as before applying the
canceller).}\label{fig_sen1}
\end{figure}

The first method proposed here for partial cancellation of NMP
zero on the Riemann surface is based on designing a canceller and
proportional controller such that the open-loop systems with and
without using canceller have the same PM while the former has a
larger DC gain. It means that application of canceller makes it
possible to use controllers with larger gains in the loop without
affecting the stability properties of system. For this purpose
consider Fig. \ref{fig_sen1} which shows the general relation
between the Bode plots of $K_{p1}P(s)$, $K_{p1}C_{canc}(s)P(s)$,
and $K_{p2}C_{canc}(s)P(s)$ for some $K_{p2}>K_{p1}$ ($K_{p1}$ and
$K_{p2}$ denote two different values for the gain of proportional
controller in Fig. \ref{fig_loop1}). Note that according to Fig.
\ref{fig_freq1} amplitude of $1-(s/z_{nmp})^\alpha$ (which appears
in the numerator of $K_{p1}C_{canc}(s)P(s)$)  is smaller than the
amplitude of $1-s/z_{nmp}$ (which appears at the numerator of
$K_{p1}P(s)$) at all frequencies. It concludes that we have
$|K_{p1} C_{canc}(j\omega) P(j\omega)|<|K_{p1} P(j\omega)|$ for
all $\omega$ as it can be observed in Fig. \ref{fig_sen1}. It also
results in the fact that the gain crossover frequency of $K_{p1}
C_{canc}(s) P(s)$ is necessarily smaller than $K_{p1}P(s)$. Recall
that all terms in the numerator and denominator of $K_{p1}
C_{canc}(s) P(s)$ and $K_{p1} P(s)$, except the term cancelled by
canceller, are exactly the same.

In order to explain the proposed design method, first assume that
in the feedback system shown in Fig. \ref{fig_loop1} we have
$C_{canc}(s)=1$ and the gain of proportional controller is chosen
such that the resulted closed-loop system has a certain PM. The
solid curve in Fig. \ref{fig_sen1} shows the Bode plot of
$K_{p1}P(s)$ where $K_{p1}$ is the gain of proportional controller
to achieve the desired PM, and points $A$ and $B$ denote the gain
crossover frequency of the resulted open-loop system and the
corresponding phase lag at this frequency, respectively. If in
this feedback system we add the canceller block in series with
plant (assuming the same value for the gain of proportional
controller and a certain value for $\alpha$) we arrive at a
feedback system whose Bode plot is shown by the dotted curve in
Fig. \ref{fig_sen1}. As mentioned earlier, the canceller has the
property that necessarily decreases the gain crossover frequency
and makes the system more phase lag at all frequencies as it can
be observed in Fig. \ref{fig_sen1}. According to this figure if we
want the feedback system with canceller has the same PM as it had
before using it, we must increase the gain of proportional
controller from $K_{p1}$ to the suitably chosen value $K_{p2}$.
More precisely, the value of $K_{p2}$ must be chosen such that the
phase lag of $K_{p2}C_{canc}(s) P(s)$ (point $C$) at its gain
crossover frequency (point $E$) be equal to the phase lag of
$K_{p1}P(s)$ (point $B$) at its gain crossover frequency (point
$A$). For this purpose, the required increment in $K_{p1}$ to
arrive at $K_{p2}$ is equal to the vertical distance between
points $D$ and $E$. In other words, after applying the canceller,
the Bode magnitude plot of the open-loop system drops and we need
to increase the gain of proportional controller to move it upward
such that the phase of open-loop system with canceller at the new
gain crossover frequency becomes equal to the one it was before
applying the canceller. Note that as it can be observed in Fig.
\ref{fig_sen1} the phase plots of $K_{p2}C_{canc}(s)P(s)$ and
$K_{p1}P(s)$ are exactly the same.

To sum up the design procedure, consider a feedback system with
the open-loop transfer function $K_{p1}P(s)$ (solid curve in Fig.
\ref{fig_sen1}) where the value of $K_{p1}$ is chosen such that
the phase lag of $K_{p1}P(s)$ (point $B$) at its gain crossover
frequency (point $A$) be equal to the desired value. Apply the
canceller in series with $K_{p1}P(s)$ assuming a certain value for
$\alpha$ to arrive at a feedback system with the open-loop
transfer function shown by dotted curve in Fig. \ref{fig_sen1}.
Then increase the gain of proportional controller from $K_{p1}$ to
$K_{p2}$ to adjust the gain crossover frequency such that the PM
becomes equal to the one it was before applying the canceller
(i.e., move point $D$ to $E$ in Fig. \ref{fig_sen1} by increasing
the gain of proportional controller from $K_{p1}$ to $K_{p2}$).
The resulted feedback system with canceller has the same PM as it
had before applying the canceller while the DC gain and the
corresponding tracking and disturbance rejection errors are
smaller. Note that using smaller values for $\alpha$ increases the
maximum of achievable DC gain at the cost of decreasing the gain
crossover frequency. Note also that one can design a controller
for the resulted augmented plant $K_{p2}C_{canc}(s)P(s)$ using any
desired method by facing less limitations caused by the NMP zero.

\begin{figure}\center
\includegraphics[width=8.5cm]{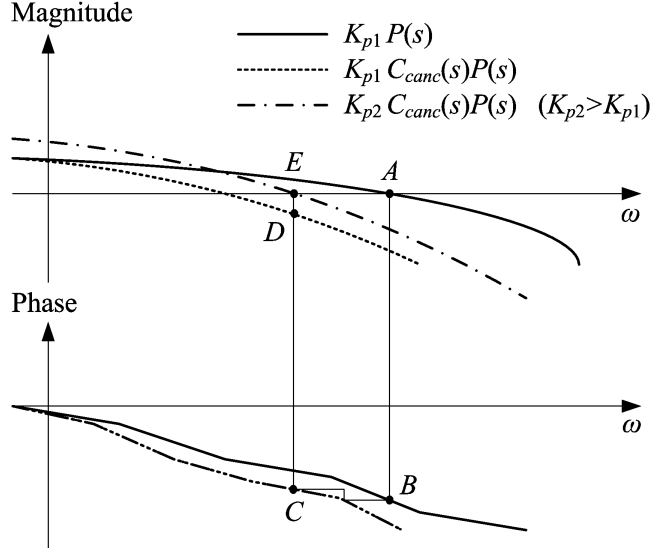}
\caption{The design procedure for canceller and proportional
controller to arrive at a feedback system with increased PM and DC
gain.}\label{fig_sen1_2}
\end{figure}

In general, it is also possible to design the proportional
controller and canceller such that both the PM and open-loop DC
gain (or equivalently, gain of proportional controller) are
increased simultaneously. The design procedure is very similar to
the previous routine and the details are shown in Fig.
\ref{fig_sen1_2}. To sum up, first determine the value of $K_{p1}$
such that $K_{p1}P(s)$ has the desired phase lag (or equivalently,
the desired PM) at the gain crossover frequency, as identified by
point $B$ in Fig. \ref{fig_sen1_2}. Next, determine point $C$ on
the Bode phase plot of $K_{p1}C_{canc}(s)P(s)$ (assuming a certain
value for $\alpha$), which is less phase lag compared to point $B$
to a desired value. More precisely, the vertical difference
between points $B$ and $C$ determines the required increment in
the PM. Then, increase the gain of proportional controller from
$K_{p1}$ to $K_{p2}$ such that point $D$ on the Bode magnitude
plot of $K_{p1}C_{canc}P(s)$ moves to point $E$ on the Bode
magnitude plot of $K_{p2}C_{canc}P(s)$. In this manner the phase
lag of $K_{p2}C_{canc}P(s)$ at its gain cross over frequency
(identified as point $E$) becomes equal to point $C$ (note that
point $C$ is less phase lag than the original system $K_{p1}P(s)$
at its gain crossover frequency). Clearly, one can also design a
controller for the resulted augmented plant $K_{p2}C_{canc}P(s)$
to arrive at a higher performance feedback system.

\subsubsection{Designing a canceller to increase PM and keep the DC gain unchanged}
\begin{figure}\center
\includegraphics[width=8.5cm]{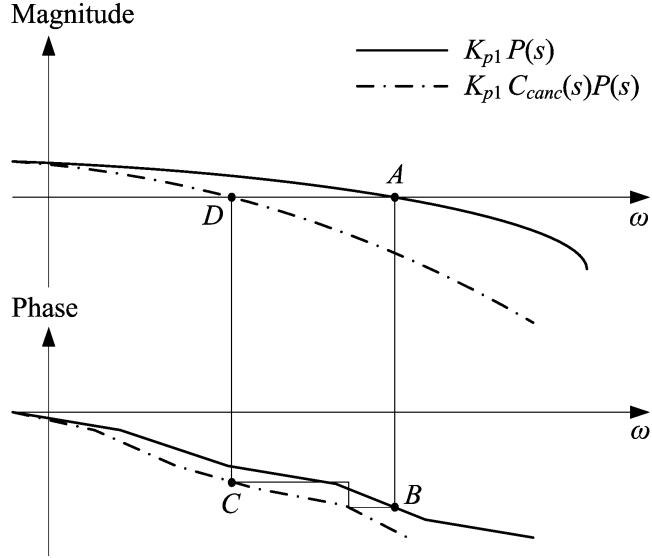}
\caption{The design procedure for canceller to increase PM and
keep the DC gain unchanged.}\label{fig_sen2}
\end{figure}

The second possible approach is to design the canceller such that
the feedback systems with and without using canceller have the
same DC gain while the system with canceller has a larger PM. For
this purpose consider again the feedback system shown in Fig.
\ref{fig_loop1} where assuming $C_{canc}(s)=1$ and $C(s)=K_{p1}$
the corresponding open-loop transfer function $K_{p1}P(s)$ has a
certain phase lag (point $B$ in Fig. \ref{fig_sen2}) at its gain
crossover frequency (point $A$ in Fig. \ref{fig_sen2}). According
to the previous discussions, application of canceller decreases
the gain crossover frequency and simultaneously makes the
open-loop system more phase lag as shown by the dash-dotted curve
in Fig. \ref{fig_sen2}. Now, consider the canceller in series with
plant (assuming the same value for the gain of proportional
controller) and determine the value of $\alpha$ by trial and error
such that the phase of $K_{p1}C_{canc}(s)P(s)$ (point $C$ in Fig.
\ref{fig_sen2}) at its gain crossover frequency (point $D$) be
larger than the phase of $K_{p1}P(s)$ (point $B$) at its gain
crossover frequency (point $A$). Clearly, the increment in PM is
equal to the vertical difference between points $B$ and $C$, which
is obtained at the cost of decreasing the closed-loop bandwidth.
Note that according to Fig. \ref{fig_sen2} the maximum possible
increment in PM is limited and strictly depends on the frequency
response of the open-loop system before applying the canceller.
Note also that after designing the canceller one can design a
suitable controller for the resulted augmented plant
$K_{p1}C_{canc}(s)P(s)$ to arrive at a desired closed-loop system.

\section{Experimental Results}
\begin{figure}\center
\includegraphics[width=8.5cm]{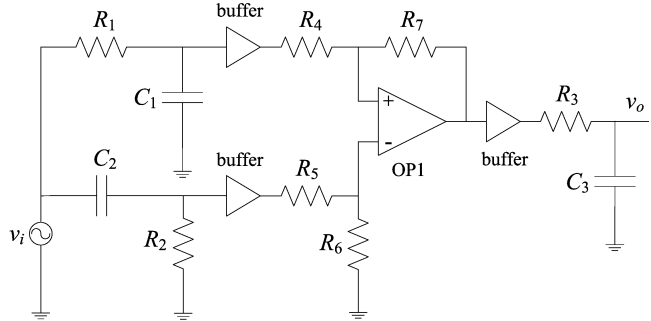}
\caption{The proposed NMP benchmark circuit with transfer function
$P(s)=V_o(s)/V_i(s)=\frac{1-R_2C_2s}{(1+R_2C_2s) (1+R_3C_3s)}$
(assuming $R_1=R_2$, $C_1=C_2$, and
$R_4=R_5=R_6=R_7$).}\label{fig_cir1}
\end{figure}

\begin{figure}\center
\includegraphics[width=5cm]{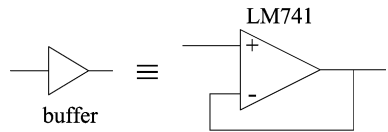}
\caption{Details of the buffers used in Fig.
\ref{fig_cir1}.}\label{fig_buffer}
\end{figure}

\begin{figure}\center
\includegraphics[width=8.5cm]{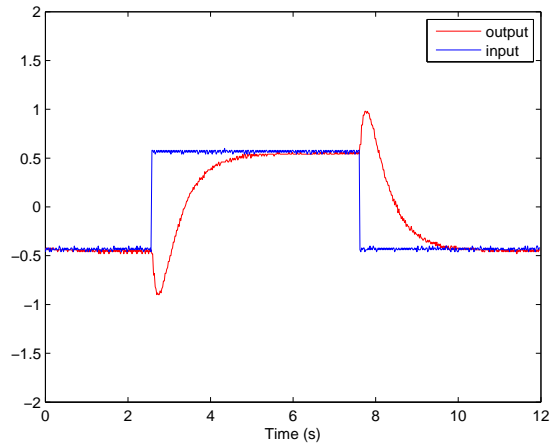}
\caption{Experimental pulse response of the proposed NMP benchmark
shown in Fig. \ref{fig_cir1}.}\label{fig_open_loop}
\end{figure}

A very linear NMP benchmark is needed for experimental
verification of the results presented in previous section. Since
most of the practical NMP systems are nonlinear to some extent,
the circuit shown in Fig. \ref{fig_cir1} is proposed in this paper
to be used as the NMP benchmark (the details of the buffers used
in this figure are shown in Fig. \ref{fig_buffer}). In this
circuit assuming $R_1=R_2$, $C_1=C_2$, and $R_4=R_5=R_6=R_7$ the
transfer function of system is calculated as the following
\begin{equation}\label{p}
P(s)=\frac{V_o(s)}{V_i(s)}=\frac{1-R_2C_2s}{(1+R_2C_2s)
(1+R_3C_3s)},
\end{equation}
which has a NMP zero at $s=1/R_2C_2$ and two stable poles at
$s=-1/R_2C_2$ and $s=-1/R_3C_3$. Note also that the DC gain of
this system is equal to unity. One advantage of this benchmark
system is that the location of its NMP zero, as well as its DC
gain and band-width, can easily be adjusted to the desired value
by assigning suitable values to resistors and capacitors. The
values of $R_1=R_2=R_4=R_5=R_6 =R_7=1.5 \mathrm{k}\Omega$,
$C_1=C_2=330\mu$F, $R_3=820\Omega$, $C_3=200\mu$F, and the
well-known op-amp LM741 are used in all simulations and the
experimental setup of this paper (the op-amps are supplied with
$\pm 10$ V DC voltage). The values assigned to resistors and
capacitors are chosen such that, firstly, the time-constant of
circuit be considerably larger than the time consumed by processor
for digital emulation of canceller, and secondly, the circuit
exhibits a considerable NMP behavior (which is identified by a
large initial undershoot in the time domain step response). More
precisely, since the canceller is realized through a very high
order FIR filter (see the discussion below) the processor should
be provided with enough time to complete the required calculations
at each sampling period. Figure \ref{fig_open_loop} shows the
experimental pulse response of the circuit shown in Fig.
\ref{fig_cir1} assuming the above mentioned values for resistors
and capacitors. As it can be observed in this figure, step
response of the open-loop system has about $46\%$ initial
undershoot which is fairly close to the one predicted by Matlab
simulation.

\begin{figure}\center
\includegraphics[width=8.5cm]{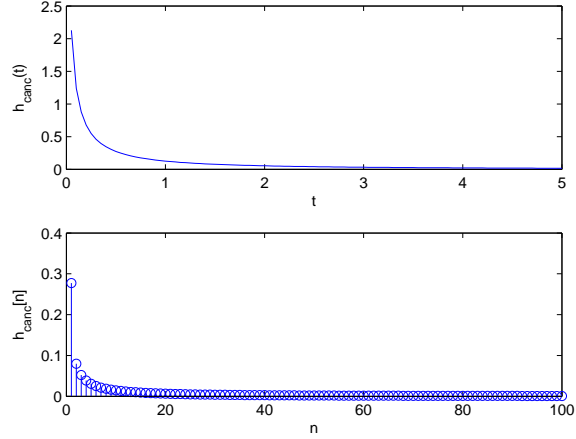}
\caption{Impulse response of canceller, $h_{canc}(t)
=L^{-1}\left\{C_{canc}(s)\right\}$, and its discrete-time
impulse-invariance equivalent, $h_{canc}[n]=T h_{canc}(nT)$
($T=50$ ms).}\label{fig_h}
\end{figure}

At this time various methods are available for realization of
simple fractional-order transfer functions like $PI^\lambda D^\mu$
or fractional-order lead-lag \citep{podlubny2} (see also \cite{farshad2}-\cite{farshad6} for more information about the simulation and tuning of fractional-order PID controllers). But, according to
the complexity of the proposed canceller these methods cannot be
used for its realization. Hence, in order to realize the canceller
first we calculate its impulse response $h_{canc}(t)$ by taking
the inverse Laplace transform from $C_{canc}(s)$, that is
\begin{equation}\label{ht}
h_{canc}(t)=L^{-1} \left\{C_{canc}(s) \right\},
\end{equation}
where according to (\ref{p}) here we have
\begin{equation}\label{canc_cir}
C_{canc}(s) = \frac{ 1}{1+(R_2C_2s)^{0.5}}.
\end{equation}
The inverse Laplace transform in (\ref{ht}) can be calculated
numerically using the Matlab function \emph{invlap.m} which can be
downloaded from Matworks website\footnote{
http://www.mathworks.com/matlabcentral/fileexchange/32824-numerical-inversion-of-laplace-transforms-in-matlab/content/INVLAP.m}.
After calculation of $h_{canc}(t)$, the equivalent discrete-time
impulse response, $h_{canc}[n]$, can be calculated using the
impulse-invariance method \citep{oppenheim} as the following
\begin{equation}\label{hn}
h_{canc}[n]=T h_{canc}(nT),
\end{equation}
where $T$ is the sampling period ($T$ is considered equal to 50 ms
in this paper). Considering the fact that the power of $s$ in
(\ref{canc_cir}) is non-integer, it is expected that $h_{canc}(t)$
decays very slowly with time \citep{podlubny}. Hence, in order to
realize $h_{canc}(t)$ with a high precision it is often needed to
approximate it with a high-order FIR filter with impulse
response $h_{canc}[n]$ as defined in (\ref{hn}). Figure
\ref{fig_h} shows $h_{canc}(t)$ and $h_{canc}[n]$ for the system under consideration where
$h_{canc}[n]$ is of length 100. Note that the DC gain of the
canceller given in (\ref{canc_cir}) is equal to unity, which
implies that the impulse response of its discrete-time equivalent
must satisfy the equality $\sum_{n=0}^\infty h_{canc}[n]=1$. But,
this condition is violated in practice since both $h_{canc}(t)$
and $h_{canc}[n]$ are necessarily truncated for realization
purposes, and moreover, the discontinuity of $h_{canc}(t)$ at
$t=0$ is the source of some errors. Hence, in order to minimize
the mismatch between theoretical and experimental results it is
better to scale all samples of the truncated discrete-time impulse
response $h_{canc}[n]$ such that their finite sum becomes equal to
unity. For this purpose all samples of the $h_{canc}[n]$
calculated from (\ref{hn}) and shown in Fig. \ref{fig_h} are
multiplied in 1.2 in the discussions of this paper. Then the
difference equation of the FIR filter with truncated impulse
response $h_{canc}[n]$ as described above is implemented using the
ATmega16 AVR microcontroller. Finally, the digital output of this
microcontroller is converted to analog using DAC08 A/D converter
and the resulted analog output is connected to the input of
benchmark circuit to form a closed-loop system. Note that in this
experiment the command signal is directly entered to the input of
microcontroller, and the subtractor and proportional controller of
the feedback system are also realized using this microcontroller.

In the following discussion whenever we talk about $P(s)$ and
$C_{canc}(s)$ equations (\ref{p}) and (\ref{canc_cir}) are under
consideration, respectively.

\begin{figure}\center
\includegraphics[width=8.5cm]{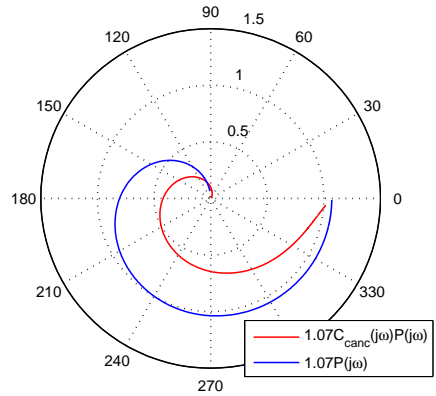}
\caption{The Nyquist plot of open-loop system with and without
using the canceller (proportional controller is applied). Both of
the open-loop systems have the same DC gain.}\label{fig_nyq1}
\end{figure}

\begin{figure}\center
\includegraphics[width=8.5cm]{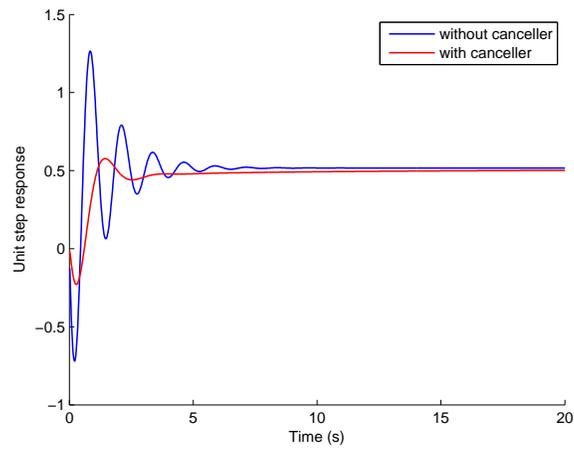}
\caption{Simulation unit step responses of the closed-loop system with and
without using canceller (a proportional controller with the
gain 1.07 is applied in both cases).}\label{fig_step1}
\end{figure}

\textbf{Scenario 1: Design for the same DC gain.} In this
experiment we design a proportional controller and canceller such
that both of the closed-loop systems (with and without using
canceller) have the same DC gain while the PM in case of using
canceller is much larger. For this purpose, without any loss of
generality, first we choose $C(s)=K_p=1.07$ and $C_{canc}(s)=1$ in
Fig. \ref{fig_loop1} to arrive at a feedback system whose PM is
approximately equal to $60^\circ$ (see Fig. \ref{fig_nyq1} for
more details). Note that since in this case the DC gain of
benchmark system is equal to unity, the DC gain of the resulted
closed-loop system is equal to $1.07/ (1+1.07)\approx 0.52$.
Similarly, by choosing $C(s)=K_p=1.07$ in Fig. \ref{fig_loop1} and
considering the canceller in the loop, the DC gain of the resulted
closed-loop system with canceller is also equal to 0.52. Figure
\ref{fig_nyq1} shows the Nyquist plot of $K_pP(s)$ and
$K_pC_{canc}(s)P(s)$ for $K_p=1.07$. As it can be observed in this
figure, application of canceller changes the frequency response of
open-loop system at higher frequencies without affecting its DC
gain. Figure \ref{fig_step1} shows the unit step response of the
closed-loop systems with and without using canceller. This figure
clearly shows that application of canceller considerably reduces
the undershoot, overshoots, and settling time of the closed-loop
step response. An explanation for the reduction of undershoot in
the step response after applying the canceller can be found in
\citep{farshad1}. The reason for considerable reduction of
overshoot in the step response after applying the canceller is
that partial cancellation highly increases the PM as it can be
observed in Fig. \ref{fig_nyq1}. More precisely, according to Fig.
\ref{fig_nyq1} the PM before and after applying the canceller is
approximately equal to $60^\circ$ and $175^\circ$, respectively.
Decreasing the settling time in Fig. \ref{fig_step1} after
applying the canceller is a direct consequent of increasing the PM
after cancellation. Note that all of the benefits observed in Fig.
\ref{fig_step1} are obtained only at the cost of slight increment
in rise time, which is a direct consequent of decreasing the gain
crossover frequency (as well as closed-loop bandwidth) after
applying the canceller.

\begin{figure}\center
\includegraphics[width=8.5cm]{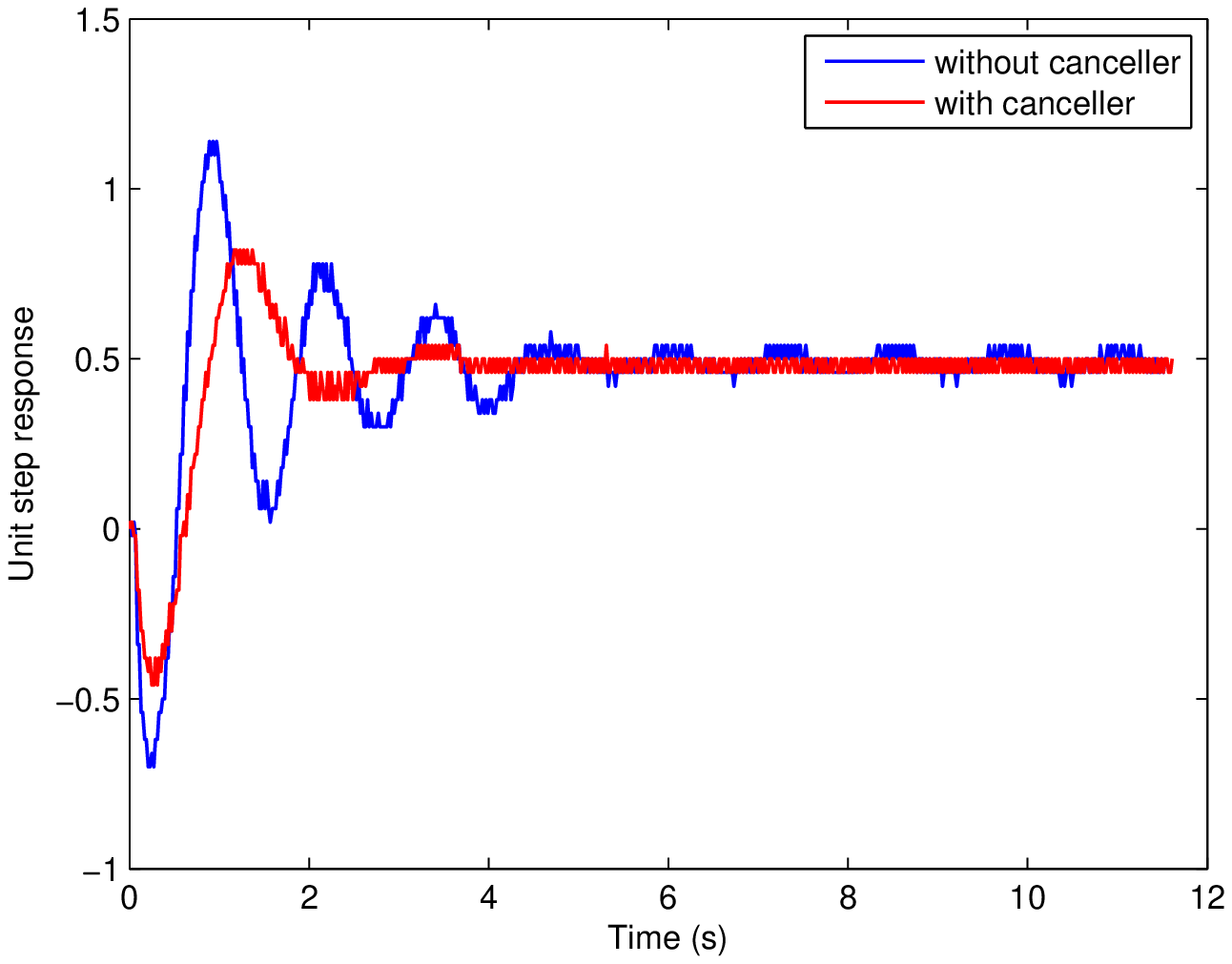}
\caption{Experimental closed-loop unit step responses with and
without using canceller (a proportional controller with the
gain 1.07 is applied in both cases).}\label{fig_exp1}
\end{figure}

The corresponding practical closed-loop step responses are shown
in Fig. \ref{fig_exp1}. The results shown in this figure are in a
fair agreement with those obtained by numerical simulations and
shown in Fig. \ref{fig_step1}. Note that there is some mismatch
between the simulation and experimental results when the canceller
is applied. This mismatch is caused because of truncation of the
impulse response of canceller as well as its discrete-time
realization.

\begin{figure}\center
\includegraphics[width=8.5cm]{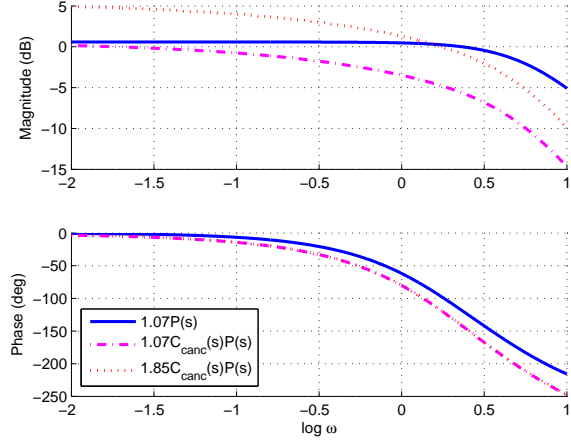}
\caption{Bode plot of $1.07P(s)$ (PM=$60^\circ$, DC gain=1.07),
$1.07 C_{canc}(s)P(s)$ (PM=$175^\circ$, DC gain=1.07), and $1.85
C_{canc}(s)P(s)$ (PM=$60^\circ$, DC
gain=1.85).}\label{fig_sen2_bode}
\end{figure}

\begin{figure}\center
\includegraphics[width=8.5cm]{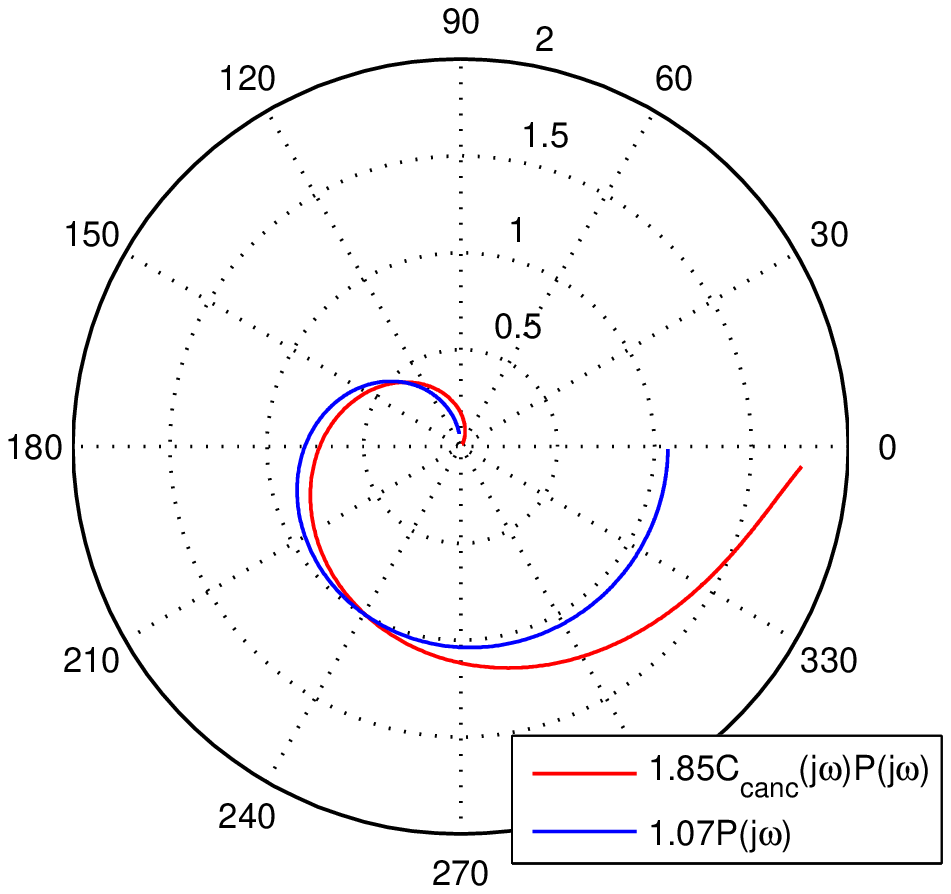}
\caption{The Nyquist plot of open-loop system with and without
using canceller (proportional controller is applied). Both of the
open-loop systems have the same PM.}\label{fig_sen2_nyq}
\end{figure}

\begin{figure}\center
\includegraphics[width=8.5cm]{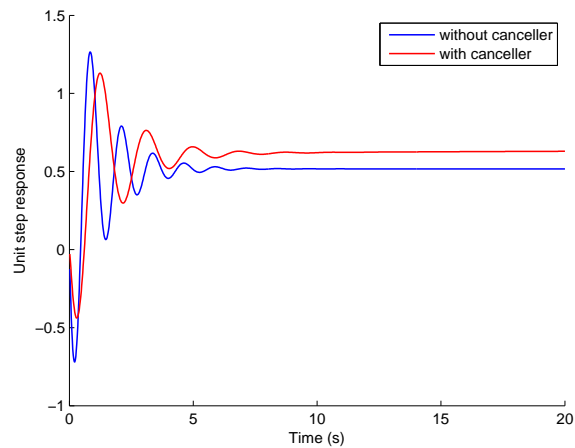}
\caption{Simulation unit step responses of the closed-loop system with and
without using canceller (the gain of proportional controller
with and without using canceller is considered equal to 1.85 and
1.07, respectively.}\label{fig_sen2_step}
\end{figure}

\begin{figure}\center
\includegraphics[width=8.5cm]{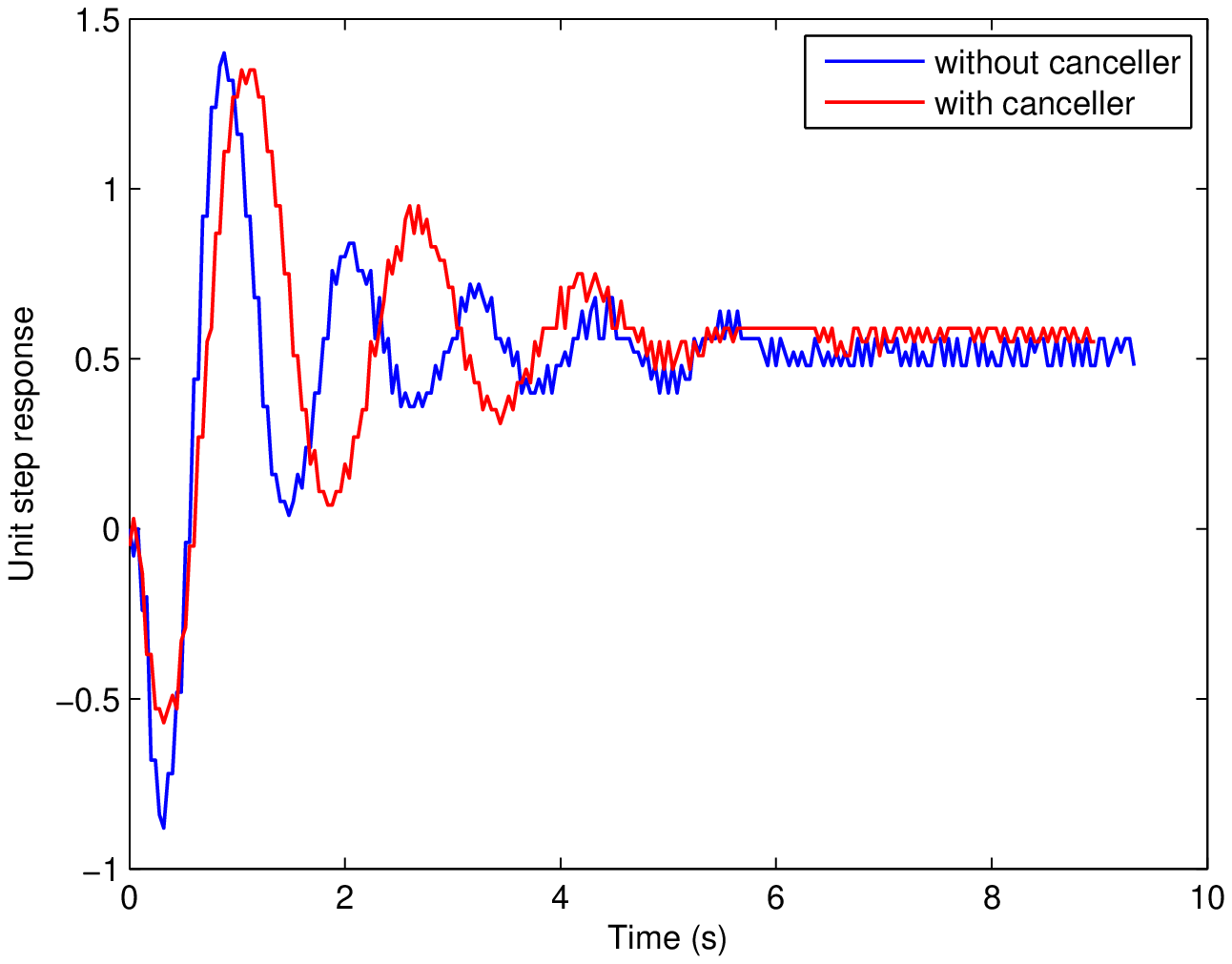}
\caption{Experimental closed-loop unit step responses with and
without using canceller (a proportional controller with the
gain 1.07 is applied in both cases).}\label{fig_exp2}
\end{figure}

\textbf{Scenario 2: Design for the same PM.} Consider again the
feedback system shown in Fig. \ref{fig_loop1} and suppose that
once without using canceller and the other time at its presence we
want to design the proportional controller such that the PM of the
resulted closed-loop systems becomes equal to $60^\circ$, and then
compare the performance of two systems. Similar to the discussions
presented in Scenario 1, when the canceller is not applied this
task can be done by choosing $K_p=1.07$. The Bode plot of the
resulted open-loop system, $1.07P(s)$, is shown by the solid curve
in Fig. \ref{fig_sen2_bode}. If in this case (i.e., assuming a
proportional controller with the gain $K_p=1.07$ in the loop) we
put the canceller given in (\ref{canc_cir}) in series with $P(s)$
we arrive at a feedback system whose open-loop Bode plot is shown
by the dash-dotted curve in Fig. \ref{fig_sen2_bode}. As it can be
observed in this figure, application of canceller without
increasing the gain of proportional controller highly increases
the PM at the cost of decreasing the bandwidth of the closed-loop
system. More precisely, application of canceller assuming
$K_p=1.07$ leads to a feedback system with $PM\approx 175^\circ$.
By a simple trial and error it can be easily verified that in case
of using canceller by choosing $K_p=1.85$ the PM of the resulted
feedback system becomes approximately equal to $60^\circ$ (see the
dotted curve in Fig. \ref{fig_sen2_bode}). In other words,
application of canceller makes it possible to use larger gains in
the loop without any reduction of PM (note that possibility of
using larger gains in the loop somehow compensates the reduction
of closed-loop bandwidth after applying canceller). Figure
\ref{fig_sen2_nyq} shows the Nyquist plots of $L_1(s) =1.07P(s)$
and $L_2(s)=1.85 C_{canc}(s) P(s)$. As it can be observed in this
figure, in both cases the PM is approximately equal to $60^\circ$
while the system with canceller has a considerably larger DC gain.
The unit step response of the closed-loop system with and without
using canceller obtained via numerical simulation is shown in Fig.
\ref{fig_sen2_step}. Note that the closed-loop system with
canceller exhibits a smaller overshoot, undershoot, and
steady-state error. The corresponding experimental closed-loop
unit step responses are shown in Fig. \ref{fig_exp2}. Although the
experimental response with canceller has some differences with the
one obtained from simulation, it still exhibits advantages
compared to the one obtained without using it. More precisely, it
can be observed in Fig. \ref{fig_exp2} that the closed-loop
response with canceller has a considerably smaller undershoot,
smaller settling time, and a smaller steady-state error compared
to the one obtained without using it. The difference between
practical and theoretical results at the presence of canceller is
firstly because of truncation of the impulse response of canceller
and secondly because of discretization which decreases the PM.

\section{Conclusion}
A method for partial cancellation of the NMP zero of a process
located in feedback connection is presented in this paper. This
partial cancellation is performed by means of a pre-compensator
(canceller) located in series with the NMP process. Two methods
for designing this pre-compensator are also proposed. In one of
these methods the pre-compensator is designed such that the PM is
increased while the open-loop DC gain remains unchanged. The other
proposed method for designing pre-compensator can simultaneously
increase the DC gain, PM, and GM of the open-loop system. Clearly,
such a change in the open-loop system by means of pre-compensator
can effectively facilitate the controller design procedure for NMP
process and makes it possible to arrive at more effective
closed-loop system.

The proposed methods for designing pre-compensator for partial
cancellation of the NMP zero of process on the Riemann surface are
also examined experimentally. For this purpose a very linear NMP
benchmark circuit is proposed and the closed-loop system
(including a proportional controller, pre-compensator, and the NMP
process) is realized using a digital micro-controller.
Experimental results show the high efficiency of the proposed
cancellation strategies.




\vspace*{6pt}
\section*{References}


\end{document}